\documentclass[leqno,draft]{article}

\def\diam{\mathop{\rm diam}}
\def\dist{\mathop{\rm dist}}

\newtheorem{theorem}{Theorem}
\newtheorem{lemma}[theorem]{Lemma}
\newtheorem{proposition}[theorem]{Proposition}
\newtheorem{definition}[theorem]{Definition}
\newtheorem{corollary}[theorem]{Corollary}

\newcommand{\begintheorem}{\addtocounter{equation}{1}\begin{theorem}}
\newcommand{\beginlemma}{\addtocounter{equation}{1}\begin{lemma}}
\newcommand{\beginproposition}{\addtocounter{equation}{1}\begin{proposition}}
\newcommand{\begindefinition}{\addtocounter{equation}{1}\begin{definition}}
\newcommand{\begincorollary}{\addtocounter{equation}{1}\begin{corollary}}

\begin{document}

\title{Elementary aspects of the geometry \\ of metric spaces}

\author{Stephen Semmes \\
        Rice University}

\date{}

\maketitle

\renewcommand{\thefootnote}{}

\footnotetext{A very basic introduction to the geometry of metric
spaces can be found in \cite{s5}, with more information in \cite{s6, s7}.}

\begin{abstract}
The setting of metric spaces is very natural for numerous questions
concerning manifolds, norms, and fractal sets, and a few of the main
ingredients are surveyed here.
\end{abstract}

\tableofcontents

\section{Metric spaces}
\label{metric spaces}
\setcounter{equation}{0}

        A \emph{metric space} is a set $M$ equipped with a function
$d(x, y)$ defined for $x, y \in M$ such that $d(x, y)$ is a
nonnegative real number that is equal to $0$ exactly when $x = y$,
\begin{equation}
        d(y, x) = d(x, y)
\end{equation}
for every $x, y \in M$, and
\begin{equation}
        d(x, z) \le d(x, y) + d(y, z)
\end{equation}
for every $x, y, z \in M$, which is known as the \emph{triangle inequality}.

        Remember that the \emph{absolute value} of a real number $r$
is denoted $|r|$ and equal to $r$ when $r \ge 0$ and to $-r$ when $r
\le 0$.  It is easy to check that
\begin{equation}
        |r + t| \le |r| + |t|
\end{equation}
and
\begin{equation}
        |r \, t| = |r| \, |t|
\end{equation}
for any pair of real numbers $r$, $t$.  The standard metric on the
real line ${\bf R}$ is given by $|r - t|$, which is the first main
example of a metric space.

        If $(M, d(x, y))$ is a metric space, then $d(x, y)$ is called
the \emph{distance function} or \emph{metric} on $M$.  For each $x \in
M$ and $r > 0$, the \emph{open ball} in $M$ with center $x$ and radius
$r$ is
\begin{equation}
        B(x, r) = \{y \in M : d(x, y) < r\}.
\end{equation}
Similarly, the \emph{closed ball} with center $x$ and radius $r \ge 0$ is
\begin{equation}
        \overline{B}(x, r) = \{y \in M : d(x, y) \le r\}.
\end{equation}
Thus
\begin{equation}
        B(x, r) \subseteq \overline{B}(x, r) \subseteq B(x, t)
\end{equation}
when $r < t$.

        Let $a$, $b$ be real numbers with $a < b$.  The \emph{open
interval} $(a, b)$ in ${\bf R}$ is defined by
\begin{equation}
        (a, b) = \{r \in {\bf R} : a < r < b\},
\end{equation}
and the \emph{closed interval} $[a, b]$ is defined by
\begin{equation}
        [a, b] = \{r \in {\bf R} : a \le r \le b\}.
\end{equation}
One may also allow $a = b$ for the latter.  The \emph{length} of these
intervals is $b - a$.  Note that open and closed balls in the real
line with respect to the standard metric are open and closed
intervals.

\section{A little calculus}
\label{calculus}
\setcounter{equation}{0}

        Suppose that $a$, $b$ are real numbers with $a < b$, and that
$f(x)$ is a continuous real-valued function on the closed interval
$[a, b]$ in the real line.  The \emph{extreme value theorem} states
that there are elements $p$, $q$ of $[a, b]$ at which $f$ attains its
maximum and minimum, which is to say that
\begin{equation}
        f(q) \le f(x) \le f(p)
\end{equation}
for every $x \in [a, b]$.  This works as well for continuous
real-valued functions on compact subsets of metric spaces, or even
topological spaces.  If $p$ or $q$ is in the open interval $(a, b)$
and $f$ is differentiable there, then the derivative $f'(p)$ or
$f'(q)$ is equal to $0$.

        Suppose that $f(x)$ is differentiable at every point in $(a,
b)$.  If $f(a) = f(b) = 0$, then \emph{Rolle's theorem} states that
$f'(x) = 0$ for some $x \in (a, b)$.  This is because the maximum or
minimum of $f$ on $[a, b]$ is attained on $(a, b)$, or $f(x) = 0$ for
every $x \in [a, b]$.  No matter the values of $f(a)$, $f(b)$, the
\emph{mean value theorem} says that there is an $x \in (a, b)$ such that
\begin{equation}
        f'(x) = \frac{f(b) - f(a)}{b - a}.
\end{equation}
This follows from Rolle's theorem applied to $f - \phi$, where
$\phi(x) = \alpha \, x + \beta$ and $\alpha, \beta \in {\bf R}$ are
chosen so that $\phi(a) = f(a)$, $\phi(b) = f(b)$.

        Of course, the derivative of a constant function is $0$, and
the mean value theorem implies that a continuous function $f$ on $[a,
b]$ is constant if the derivative of $f$ exists and is equal to $0$ at
every point in $(a, b)$.  If $f$ is monotone increasing on $[a, b]$,
in the sense that $f(x) \le f(y)$ when $a \le x \le y \le b$, then
$f'(x) \ge 0$ for every $x \in (a, b)$ at which $f$ is differentiable.
Conversely, if $f$ is continuous on $[a, b]$, differentiable on $(a,
b)$, and $f'(x) \ge 0$ for each $x \in (a, b)$, then $f$ is monotone
increasing on $[a, b]$, by the mean value theorem.  If $f'(x) > 0$ for
every $x \in (a, b)$, then $f$ is strictly increasing on $[a, b]$, in
the sense that $f(w) < f(y)$ when $a \le w < y \le b$.  However, the
derivative of a strictly increasing function may be equal to $0$, as
when $f(x) = x^3$.

\section{Norms on ${\bf R}^n$}
\label{norms}
\setcounter{equation}{0}

        Let $n$ be a positive integer, and let ${\bf R}^n$ be the
space of $n$-tuples of real numbers.  This means that an element $x$
of ${\bf R}^n$ is of the form $x = (x_1, \ldots, x_n)$, where the
coordinates $x_1, \ldots, x_n$ of $x$ are real numbers.  Addition and
scalar multiplication on ${\bf R}^n$ are defined coordinatewise in the
usual way, so that ${\bf R}^n$ becomes a finite-dimensional vector
space over the real numbers.

        A \emph{norm} on ${\bf R}^n$ is a function $N(x)$ such that
$N(x)$ is a nonnegative real number for every $x \in {\bf R}^n$ which
is equal to $0$ exactly when $x = 0$,
\begin{equation}
\label{N(r x) = |r| N(x)}
        N(r \, x) = |r| \, N(x)
\end{equation}
for every $r \in {\bf R}$ and $x \in {\bf R}^n$, and
\begin{equation}
\label{N(x + y) le N(x) + N(y)}
        N(x + y) \le N(x) + N(y)
\end{equation}
for every $x, y \in {\bf R}^n$.  If $N$ is a norm on ${\bf R}^n$, then
\begin{equation}
        d_N(x, y) = N(x - y)
\end{equation}
is a metric on ${\bf R}^n$.

        For example, the absolute value function is a norm on ${\bf
R}$, for which the corresponding metric is the standard metric on the
real line.  The standard Euclidean norm on ${\bf R}^n$ is defined by
\begin{equation}
        |x| = \Big(\sum_{j = 1}^n x_j^2 \Big)^{1/2},
\end{equation}
and the corresponding metric is the standard Euclidean metric on ${\bf R}^n$.
It is not so obvious that this satisfies the triangle inequality, and hence
is a norm, and we shall discuss a proof of this fact in Section \ref{convex 
sets}.

        One can check directly that
\begin{equation}
        \|x\|_1 = \sum_{j = 1}^n |x_j|
\end{equation}
and
\begin{equation}
        \|x\|_\infty = \max(|x_1|, \ldots, |x_n|)
\end{equation}
are norms on ${\bf R}^n$.  We shall see in Section \ref{convex sets} that
\begin{equation}
        \|x\|_p = \Big(\sum_{j = 1}^n |x_j|^p \Big)^{1/p}
\end{equation}
is a norm when $p \ge 1$, which includes the Euclidean norm as a
special case.

\section{Convex functions}
\label{convex functions}
\setcounter{equation}{0}

        A real-valued function $f(x)$ on the real line is said to be
\emph{convex} if
\begin{equation}
\label{convexity}
        f(t \, x + (1 - t) \, y) \le t \, f(x) + (1 - t) \, f(y)
\end{equation}
for every $x, y \in {\bf R}$ and $t \in [0, 1]$.  This is equivalent
to
\begin{equation}
\label{convexity, 2}
        \frac{f(w) - f(x)}{w - x} \le \frac{f(y) - f(w)}{y - w}
\end{equation}
for every $x, w, y \in {\bf R}$ such that $x < w < y$.  Applying this
condition twice, we get that
\begin{equation}
\label{convexity, 3}
        \frac{f(w) - f(x)}{w - x} \le \frac{f(z) - f(y)}{z - y}
\end{equation}
when $x < w < y < z$.  As another refinement of (\ref{convexity, 2}),
one can use (\ref{convexity}) to show that
\begin{equation}
\label{convexity, 4}
        \frac{f(w) - f(x)}{w - x} \le \frac{f(y) - f(x)}{y - x}
                                   \le \frac{f(y) - f(w)}{y - w}
\end{equation}
when $x < w < y$.

        If $f$ is differentiable and $f'$ is monotone increasing, then
the mean value theorem implies (\ref{convexity, 2}) and hence that $f$
is convex.  Conversely, (\ref{convexity, 3}) implies that the
derivative of $f$ is monotone increasing when $f$ is differentiable.
Actually, one can show that the right and left derivatives $f_+'(x)$,
$f_-'(x)$ exist for each $x \in {\bf R}$ when $f$ is convex, and
satisfy
\begin{equation}
        f_-'(x) \le f_+'(x)
\end{equation}
and
\begin{equation}
        f_+'(x) \le f_-'(y)
\end{equation}
when $x < y$.  One can also show that these conditions characterize
convexity, using analogues of Rolle's theorem and the mean value
theorem for functions with one-sided derivatives.

        A function $f : {\bf R} \to {\bf R}$ is \emph{strictly convex} if
\begin{equation}
        f(t \, x + (1 - t) \, y) < t \, f(x) + (1 - t) \, f(y)
\end{equation}
when $x \ne y$ and $0 < t < 1$.  This corresponds to strict inequality
in (\ref{convexity, 2}), (\ref{convexity, 3}), and (\ref{convexity,
4}) as well.  If $f$ is differentiable on ${\bf R}$, then $f$ is
strictly convex if and only if $f'$ is strictly increasing.
Otherwise, strict convexity can be characterized in terms of one-sided
derivatives by the requirement that
\begin{equation}
        f_+'(x) < f_-'(y)
\end{equation}
when $x < y$.  Alternatively, if a convex function $f$ on ${\bf R}$ is
not strictly convex, then $f$ is equal to an affine function on an
interval of positive length.

        For example, consider $f(r) = |r|^p$, $p > 0$.  If $p = 1$,
then $f(r) = |r|$ is convex but not strictly convex on ${\bf R}$.  If
$p = 2$, then $f(r) = r^2$ is twice-differentiable, $f''(r) = 2$, and
$f$ is strictly convex.  If $p > 2$, then $f$ is twice-differentiable
on ${\bf R}$, $f''(r) > 0$ when $r \ne 0$, $f''(0) = 0$, and $f$ is
strictly convex because $f'$ is strictly increasing.  If $1 < p < 2$,
then $f$ is differentiable on ${\bf R}$, twice-differentiable on ${\bf
R} \backslash \{0\}$, $f''(r) > 0$ when $r \ne 0$, and again $f$ is
strictly convex since $f'$ is strictly increasing.  If $0 < p < 1$,
then $f$ is twice-differentiable on ${\bf R} \backslash \{0\}$,
$f''(r) < 0$ when $r \ne 0$, and $f$ is not convex.

\section{Convex sets}
\label{convex sets}
\setcounter{equation}{0}

        A set $E \subseteq {\bf R}^n$ is said to be \emph{convex} if
\begin{equation}
        t \, x + (1 - t) \, y \in E
\end{equation}
for every $x, y \in E$ and $t \in (0, 1)$.  For example, open and
closed balls associated to metrics defined by norms on ${\bf R}^n$ are
convex.

        Conversely, suppose that $N(x)$ is a nonnegative real-valued
function on ${\bf R}^n$ such that $N(x) > 0$ when $x \ne 0$ and the
homogeneity condition (\ref{N(r x) = |r| N(x)}) holds for all $x \in
{\bf R}^n$ and $r \in {\bf R}$.  If the closed unit ball
\begin{equation}
\label{B_N}
        B_N = \{x \in {\bf R}^n : N(x) \le 1\}
\end{equation}
is convex, then $N$ satisfies the triangle inequality (\ref{N(x + y)
le N(x) + N(y)}) and hence is a norm.  Let $x, y \in {\bf R}^n$ be
given, and let us check (\ref{N(x + y) le N(x) + N(y)}).  We may
suppose that $x, y \ne 0$, since the inequality is trivial when $x =
0$ or $y = 0$.  Put
\begin{equation}
        x' = \frac{x}{N(x)}, \quad y' = \frac{y}{N(y)},
\end{equation}
so that $N(x') = N(y') = 1$.  By hypothesis,
\begin{equation}
        N(t \, x' + (1 - t) \, y') \le 1
\end{equation}
when $0 \le t \le 1$.  Applying this with
\begin{equation}
        t = \frac{N(x)}{N(x) + N(y)},
\end{equation}
we get (\ref{N(x + y) le N(x) + N(y)}), as desired.

        For example, suppose that $N(x) = \|x\|_p$, $1 < p < \infty$.
Let $x, y \in {\bf R}^n$ with $\|x\|_p, \|y\|_p \le 1$ be given, so that
\begin{equation}
        \sum_{j = 1}^n |x_j|^p, \ \sum_{j = 1}^n |y_j|^p \le 1.
\end{equation}
We would like to show that
\begin{equation}
        \|t \, x + (1 - t) \, y\|_p \le 1
\end{equation}
when $0 \le t \le 1$, which is the same as
\begin{equation}
        \sum_{j = 1}^n |t \, x_j + (1 - t) \, y_j|^p \le 1.
\end{equation}
The convexity of $|r|^p$ on ${\bf R}$ implies that
\begin{equation}
        |t \, x_j + (1 - t) \, y_j|^p \le t |x_j|^p + (1 - t) \, |y_j|^p
\end{equation}
for each $j$, and the desired inequality follows by summing this over
$j$.

        A norm $N$ on ${\bf R}^n$ is said to be \emph{strictly convex}
if the unit ball $B_N$ is strictly convex in the sense that
\begin{equation}
        N(t \, x + (1 - t) \, y) < 1
\end{equation}
when $x, y \in {\bf R}^n$, $N(x) = N(y) = 1$, $x \ne y$, and $0 < t <
1$.  It is easy to see that the absolute value function is strictly
convex as a norm on ${\bf R}$, if not as a general function as in the
previous section.  One can also check that $\|x\|_p$ is a strictly
convex norm on ${\bf R}^n$ when $p > 1$, using the strict convexity of
$|r|^p$ on ${\bf R}$ and computations as in the preceding paragraph.
However, $\|x\|_1$ and $\|x\|_\infty$ are not strictly convex norms on
${\bf R}^n$ when $n \ge 2$.

\section{A little more calculus}
\label{a little more calculus}
\setcounter{equation}{0}

        Let $f$ be a continuous real-valued function on a closed
interval $[a, b]$, $a < b$.  The integral
\begin{equation}
        \int_a^b f(t) \, dt
\end{equation}
can be defined in the usual way as a limit of finite sums.  The
convergence of the finite sums to the integral uses the fact that
continuous functions on $[a, b]$ are actually uniformly continuous.
It is well known that continuous functions on compact subsets of any
metric space are uniformly continuous.

        Consider the indefinite integral
\begin{equation}
        F(x) = \int_a^x f(t) \, dt.
\end{equation}
This defines a continuous function on $[a, b]$ which is differentiable
on $(a, b)$ and satisfies $F'(x) = f(x)$.  Similarly, $F$ has
one-sided derivatives at the endpoints $a$, $b$ that satisfy the same
condition.  If another differentiable function on $(a, b)$ has
derivative $f$, then the difference of $F$ and this function is
constant, by the mean value theorem.

        Clearly $F$ is monotone increasing on $[a, b]$ if $f \ge 0$ on
the whole interval.  If $f > 0$ on $[a, b]$, then $F$ is strictly
increasing.  The same conclusion holds if $f \ge 0$ on $[a, b]$ and $f
> 0$ at some point in any nontrivial subinterval of $[a, b]$.
Equivalently, if $f \ge 0$ on $[a, b]$, and if $F$ is not strictly
increasing on $[a, b]$, then $f = 0$ at every point in a nontrivial
subinterval.

        Suppose that $f' \ge 0$ on $(a, b)$, or simply that $f$ is
monotone increasing on $[a, b]$.  This implies that
\begin{equation}
        f(x) \le \frac{F(y) - F(x)}{y - x} \le f(y)
\end{equation}
when $a \le x < y \le b$.  In particular,
\begin{equation}
        \frac{F(w) - F(x)}{w - x} \le \frac{F(y) - F(w)}{y - w}
\end{equation}
when $a \le x < w < y \le b$.  If $f$ is strictly increasing, then
these inequalities are strict as well.

\section{Supremum and infimum}
\label{supremum and infimum}
\setcounter{equation}{0}

        A real number $b$ is said to be an \emph{upper bound} for a
set $A \subseteq {\bf R}$ if $a \le b$ for every $a \in A$.  We say
that $b_1 \in {\bf R}$ is the \emph{least upper bound} or
\emph{supremum} of $A$ if $b_1$ is an upper bound for $A$ and $b_1 \le
b$ for every upper bound $b$ of $A$.  If $b_2 \in {\bf R}$ also
satisfies these two conditions, then $b_1 \le b_2$ and $b_2 \le b_1$,
and hence $b_1 = b_2$.  Thus the supremum of $A$ is unique when it
exists, in which case it is denoted $\sup A$.  The \emph{completeness
property} of the real line states that every nonempty set with an
upper bound has a least upper bound.

        More precisely, this is completeness with respect to the
ordering on the real line, which can be defined for other ordered
sets.  There is also completeness for metric spaces, which means that
every Cauchy sequence converges.  Both forms of completeness hold on
the real line, and are basically equivalent to each other in this
particular situation.  However, the two notions are distinct, because
they can be applied in different circumstances.  There are
completeness conditions concerning the existence of solutions of
ordinary differential equations as well, which may be related to
completeness for an associated metric space.

       Similarly, a real number $c$ is said to be a \emph{lower bound}
for $A \subseteq {\bf R}$ if $c \le a$ for every $a \in A$, and $c_1
\in {\bf R}$ is the \emph{greatest lower bound} or \emph{infimum} of
$A$ if $c_1$ is a lower bound for $A$ and $c \le c_1$ for every lower
bound $c$ of $A$.  This is unique when it exists for the same reasons
as before, and is denoted $\inf A$.  It follows from completeness that
a nonempty set $A \subseteq {\bf R}$ with a lower bound has a greatest
lower bound, which can be characterized as the supremum of the set of
lower bounds of $A$.  Alternatively, the infimum of $A$ is equal to
the negative of the supremum of $-A = \{-a : a \in A\}$.

\section{Bounded sets}
\label{bounded sets}
\setcounter{equation}{0}

        Let $(M, d(x, y))$ be a metric space.  A set $E \subseteq M$
is said to be \emph{bounded} if there is a $p \in M$ and an $r \ge 0$
such that
\begin{equation}
        d(p, x) \le r
\end{equation}
for every $x \in E$.  This implies that for every $q \in M$ there is a
$t \ge 0$ such that $d(q, x) \le t$ for every $x \in E$, by taking $t
= r + d(p, q)$.

        Equivalently, $E \subseteq M$ is bounded if the set of
distances $d(x, y)$ for $x, y \in E$ has an upper bound in ${\bf R}$.
If $E$ is nonempty and bounded, then the \emph{diameter} of $E$ is defined by
\begin{equation}
        \diam E = \sup \{d(x, y) : x, y \in E\}.
\end{equation}
The diameter of the empty set may be interpreted as $0$.

        If $E_1 \subseteq E_2 \subseteq M$ and $E_2$ is bounded, then
$E_1$ is bounded, and
\begin{equation}
        \diam E_1 \le \diam E_2.
\end{equation}
The union of two bounded subsets of $M$ is also bounded, but the diameter of
the union may be much larger than the sum of the diameters of the two subsets.

        Suppose that $M$ is ${\bf R}^n$ equipped with a norm $N$ and
its associated metric $d_N(x, y)$.  The \emph{convex hull}
$\widehat{E}$ of a set $E \subseteq {\bf R}^n$ consists of all convex
combinations of elements of $E$.  More precisely, $\widehat{E}$ is the
set of all finite sums of the form
\begin{equation}
        \sum_{i = 1}^k r_i \, x(i),
\end{equation}
where $k$ is a positive integer, $r_1, \ldots, r_k$ are
nonnegative real numbers such that
\begin{equation}
        \sum_{i = 1}^k r_k = 1,
\end{equation}
and $x(1), \ldots, x(k)$ are elements of $E$.  It is well known that
one can take $k = n + 1$ here, but we shall not need this fact.  By
construction, $\widehat{E}$ is a convex set in ${\bf R}^n$ that
contains $E$.  Moreover, $\widehat{E}$ is the smallest such set, in
the sense that $\widehat{E}$ is contained in any convex set in ${\bf
R}^n$ that contains.  If $E$ is bounded, so that $E$ is contained in a
ball, then $\widehat{E}$ is contained in the same ball, and hence
$\widehat{E}$ is bounded.  Let us check that
\begin{equation}
        {\diam}_N \widehat{E} \le {\diam}_N E,
\end{equation}
where the subscript $N$ indicates that the diameter uses the norm $N$.
Let
\begin{equation}
 \xi = \sum_{i = 1}^k r_i \, x(i), \quad \eta = \sum_{j = 1}^l t_j \, y(j)
\end{equation}
be arbitrary elements of $\widehat{E}$, as before.  Thus
\begin{equation}
        \xi - \eta = \sum_{i = 1}^k \sum_{j = 1}^l r_i \, t_j \, (x(i) - y(j)),
\end{equation}
and therefore
\begin{equation}
 N(\xi - \eta) \le \sum_{i = 1}^k \sum_{j = 1}^l r_i \, t_j \, N(x(i) - y(j)),
\end{equation}
by the properties of norms.  This implies that
\begin{equation}
 N(\xi - \eta) \le \max \{N(x(i) - y(j)) : 1 \le i \le k, \, 1 \le j \le l\},
\end{equation}
and consequently $N(\xi - \eta) \le {\diam}_N E$, as desired.

\section{Lipschitz mappings}
\label{lipschitz mappings}
\setcounter{equation}{0}

        Let $(M_1, d_1(x, y))$ and $(M_2, d_2(u, v))$ be metric
spaces.  A mapping $f : M_1 \to M_2$ is said to be \emph{Lipschitz} if
\begin{equation}
        d_2(f(x), f(y)) \le C \, d_1(x, y)
\end{equation}
for some $C \ge 0$ and all $x, y \in M$.  More precisely, this means
that $f$ is Lipschitz of order $1$, and we shall discuss other
Lipschitz conditions later.  One can also say that $f$ is
$C$-Lipschitz or $C$-Lipschitz of order $1$ to mention the constant
$C$ explicitly.

         Thus $f$ is $C$-Lipschitz with $C = 0$ if and only if $f$ is
constant.  Note that Lipschitz mappings are uniformly continuous.
Suppose that $(M_3, d_3(w, z))$ is another metric space, and that $f_1
: M_1 \to M_2$ and $f_2 : M_2 \to M_3$ are Lipschitz mappings with
constants $C_1$, $C_2$, respectively.  The composition $f_2 \circ f_1$
is the mapping from $M_1$ to $M_2$ defined by
\begin{equation}
         (f_2 \circ f_1)(x) = f_2(f_1(x)),
\end{equation}
and it is easy to check that this is Lipschitz with constant equal to
the product of $C_1$ and $C_2$.

         If $f : M_1 \to M_2$ is $C$-Lipschitz and $E \subseteq M_1$
is bounded, then
\begin{equation}
         f(E) = \{f(x) : x \in E\}
\end{equation}
is bounded in $M_2$, and
\begin{equation}
         {\diam}_2 f(E) \le C \, {\diam}_1 E.
\end{equation}
Here the subscripts indicate in which metric space the diameter is
taken.  This is easy to verify, directly from the definitions, and
suggests another way to look at the composition of Lipschitz mappings,
as in the previous paragraph.

\section{Real-valued functions}
\label{real functions}
\setcounter{equation}{0}

        Let $f$ be a real-valued function on an open interval $(a, b)$
in the real line.  If $f$ is $C$-Lipschitz with respect to the
standard metric on the domain and range, then
\begin{equation}
        |f'(x)| \le C
\end{equation}
at every point $x \in (a, b)$ at which $f$ is differentiable, by
definition of the derivative.  Conversely, if $f$ is differentiable
and satisfies this condition everywhere on $(a, b)$, then $f$ is
$C$-Lipschitz, by the mean value theorem.

        Now let $(M, d(x, y))$ be a metric space.  A function $f : M
\to {\bf R}$ is $C$-Lipschitz with respect to the standard metric on
${\bf R}$ if and only if
\begin{equation}
        f(x) \le f(y) + C \, d(x, y)
\end{equation}
for every $x, y \in M$.  This follows easily from the definitions.  In
particular, $f_p(x) = d(p, x)$ is $1$-Lipschitz for every $p \in M$.

        If $A \subseteq M$, $A \ne \emptyset$, and $x \in M$, then put
\begin{equation}
        \dist(x, A) = \inf \{d(x, a) : a \in A\}.
\end{equation}
For each $x, y \in M$ and $a \in A$,
\begin{equation}
        \dist(x, A) \le d(x, a) \le d(x, y) + d(y, a),
\end{equation}
and therefore
\begin{equation}
        \dist(x, A) \le \dist(y, A) + d(x, y).
\end{equation}
This shows that $\dist(x, A)$ is $1$-Lipschitz on $M$.

         Suppose that $f_1, f_2 : M \to {\bf R}$ are Lipschitz with
constants $C_1$, $C_2$, respectively.  For any $r_1, r_2 \in {\bf R}$,
$r_1 \, f_1 + r_2 \, f_2$ is Lipschitz with constant $|r_1| \, C_1 +
|r_2| \, C_2$.  Suppose also that $f_1$, $f_2$ are bounded on $M$, with
\begin{equation}
        |f_1(x)| \le k_1, \quad |f_2(x)| \le k_2
\end{equation}
for some $k_1, k_2 \ge 0$ and every $x \in M$.  Because
\begin{eqnarray}
\lefteqn{f_1(x) \, f_2(x) - f_1(y) \, f_2(y)}   \\
 & & = (f_1(x) - f_1(y)) \, f_2(x) + f_1(y) \, (f_2(x) - f_2(y)) \nonumber
\end{eqnarray}
for every $x, y \in M$, $f_1 \, f_2$ is Lipschitz on $M$ with constant
$k_2 \, C_1 + k_1 \, C_2$.

\section{${\bf R}^n$-valued functions}
\label{R^n functions}
\setcounter{equation}{0}

        Let $N$ be a norm on ${\bf R}^n$.  Thus $N$ is $1$-Lipschitz
as a real-valued function on ${\bf R}^n$ with the metric $d_N(x, y)$
associated to $N$, as in the previous section.  One can also show that
$N$ is bounded by a constant multiple of the standard Euclidean norm
on ${\bf R}^n$.  This uses the finite-dimensionality of ${\bf R}^n$ in
an essential way, and it implies that $N$ is Lipschitz with respect to
the standard metric on ${\bf R}^n$.

        Suppose that $f$ is a continuous ${\bf R}^n$-valued function
on a closed interval $[a, b]$ in the real line.  As an extension of
the triangle inequality for $N$,
\begin{equation}
        N\Big(\int_a^b f(t) \, dt\Big) \le \int_a^b N(f(t)) \, dt.
\end{equation}
Indeed, the analogous statement for the finite sums follows from the
triangle inequality for $N$.  The integral of $f$ can be approximated
by finite sums, and continuity of $N$ as in the preceding paragraph
can be employed to pass to the limit.  Alternatively, one can use
duality, as follows.  For any linear functional $\phi : {\bf R}^n \to
{\bf R}$,
\begin{equation}
        \phi\Big(\int_a^b f(t) \, dt\Big) = \int_a^b \phi(f(t)) \, dt.
\end{equation}
If $|\phi(w)| \le N(w)$ for every $w \in {\bf R}^n$, then we get that
\begin{equation}
        \biggl|\int_a^b \phi(f(t)) \, dt\biggr| \le \int_a^b N(f(t)) \, dt.
\end{equation}
A famous theorem states that for each $v \in {\bf R}^n$ there is such
a $\phi$ with $\phi(v) = N(v)$, which permits one to estimate the norm
of the integral.  We shall not discuss the proof of this here, but one
can take $\phi(w)$ to be the standard inner product of $w$ with $v /
|v|$ when $v \ne 0$ and $N$ is the Euclidean norm on ${\bf R}^n$, and
there are also explicit expressions for $\phi$ when $N(w) = \|w\|_p$,
$1 \le p \le \infty$.

        Suppose now that $F : [a, b] \to {\bf R}^n$ is $C$-Lipschitz
with respect to the standard metric on ${\bf R}$ and the metric $d_N$
on ${\bf R}^n$.  If $F$ is differentiable at a point $x \in (a, b)$,
then $N(F'(x)) \le C$.  This follows from the definition of the
derivative, as in the real-valued case.  Conversely, if $F$ is
continuously differentiable on $[a, b]$ and $N(F') \le C$, then one
can use the fundamental theorem of calculus and the integral form of
the triangle inequality to show that that $F$ is $C$-Lipschitz with
respect to $N$.  One can use duality to get the same conclusion when
$F$ is continuous on $[a, b]$ and differentiable on $(a, b)$ with
$N(F') \le C$, by applying the mean value theorem to $\phi \circ F$
for linear functionals $\phi : {\bf R}^n \to {\bf R}$.

        Let $(M, d(x, y))$ be a metric space, and let $F = (F_1,
\ldots, F_n)$ be a mapping from $M$ into ${\bf R}^n$.  If ${\bf R}^n$
is equipped with the norm $\|w\|_\infty$, then it is easy to see that
$F$ is $C$-Lipschitz if and only if $F_1, \ldots, F_n$ are
$C$-Lipschitz as real-valued functions on $M$.  Of course, one can
estimate Lipschitz conditions for $F$ in terms of Lipschitz conditions
for $F_1, \ldots, F_n$ for other norms on ${\bf R}^n$, and vice-versa,
but the relationship between the constants is normally not quite as
simple as for the norm $\|w\|_\infty$.

\section{Bounded variation}
\label{bounded variation}
\setcounter{equation}{0}

        Let $f$ be a real-valued function on a closed interval $[a,
b]$.  A \emph{partition} of $[a, b]$ is a finite sequence $\{t_j\}_{j
= 0}^n$ of real numbers such that
\begin{equation}
        a = t_0 < t_1 < \cdots < t_n = b.
\end{equation}
For each partition $\mathcal{P} = \{t_j\}_{j = 0}^n$ of $[a, b]$, consider
\begin{equation}
        V_\mathcal{P}(f) = \sum_{j = 1}^n |f(t_j) - f(t_{j - 1})|.
\end{equation}
This measures the variation of $f$ on the partition $\mathcal{P}$.  We
say that $f$ has \emph{bounded variation} on $[a, b]$ if there is an
upper bound for $V_\mathcal{P}(f)$ over all partitions $\mathcal{P}$
of $[a, b]$.  In this case, the \emph{total variation} $V_a^b(f)$ of
$f$ on $[a, b]$ is defined by
\begin{equation}
\label{V_a^b(f)}
        V_a^b(f) = \sup \{ V_\mathcal{P} : \mathcal{P}
                                    \hbox{ is a partition of } [a, b]\}.
\end{equation}
Thus $V_a^b(f) = 0$ if and only if $f$ is constant on $[a, b]$.

        Using the partition that consists of only $a$, $b$, we get that
\begin{equation}
        |f(b) - f(a)| \le V_a^b(f).
\end{equation}
If $f$ is monotone increasing on $[a, b]$, then
\begin{equation}
        V_\mathcal{P}(f) = f(b) - f(a)
\end{equation}
for every partition $\mathcal{P}$ of $[a, b]$.  Hence $f$ has bounded
variation on $[a, b]$, and
\begin{equation}
        V_a^b(f) = f(b) - f(a).
\end{equation}
Conversely, if $f$ has bounded variation on $[a, b]$ and
\begin{equation}
        V_a^b(f) = |f(b) - f(a)|,
\end{equation}
then $f$ is either monotone increasing or decreasing on $[a, b]$.

        If $f$ is $C$-Lipschitz on $[a, b]$, then
\begin{equation}
        V_\mathcal{P}(f) \le C \, (b - a)
\end{equation}
for every partition $\mathcal{P}$ of $[a, b]$.  Hence $f$ has bounded
variation on $[a, b]$, and
\begin{equation}
        V_a^b(f) \le C \, (b - a).
\end{equation}
If $\phi : {\bf R} \to {\bf R}$ is $C$-Lipschitz, then
\begin{equation}
        V_\mathcal{P}(\phi \circ f) \le C \, V_\mathcal{P}(f)
\end{equation}
for every $f : [a, b] \to {\bf R}$ and partition $\mathcal{P}$ of $[a, b]$.
If $f$ has bounded variation on $[a, b]$, then it follows that $\phi
\circ f$ has bounded variation on $[a, b]$, and
\begin{equation}
        V_a^b(\phi \circ f) \le C \, V_a^b(f).
\end{equation}

        Let $f_1, f_2 : [a, b] \to {\bf R}$ and $r_1, r_2 \in {\bf R}$
be given.  For any partition $\mathcal{P}$ of $[a, b]$,
\begin{equation}
        V_\mathcal{P}(r_1 \, f_1 + r_2 \, f_2)
           \le |r_1| \, V_\mathcal{P}(f_1) + |r_2| \, V_\mathcal{P}(f_2).
\end{equation}
If $f_1$, $f_2$ have bounded variation on $[a, b]$, then it follows
that $r_1 \, f_1 + r_2 \, f_2$ also has bounded variation, with
\begin{equation}
 V_a^b(r_1 \, f_1 + r_2 \, f_2) \le |r_1| \, V_a^b(f_1) + |r_2| \, V_a^b(f_2).
\end{equation}

        Suppose that $f_1$, $f_2$ are bounded on $[a, b]$, so that
\begin{equation}
        |f_1(x)| \le k_1, \quad |f_2(x)| \le k_2
\end{equation}
for some $k_1, k_2 \ge 0$ and every $x \in [a, b]$.  It is easy to
check that
\begin{equation}
        V_\mathcal{P}(f_1 \, f_2) \le k_2 \, V_\mathcal{P}(f_1)
                                     + k_1 \, V_\mathcal{P}(f_2)
\end{equation}
for every partition $\mathcal{P}$ of $[a, b]$.  If $f_1$, $f_2$ have
bounded variation on $[a, b]$, then $f_1 \, f_2$ has bounded variation, and
\begin{equation}
        V_a^b(f_1 \, f_2) \le k_2 \, V_a^b(f_1) + k_1 \, V_a^b(f_2).
\end{equation}
This is analogous to the earlier estimate for the Lipschitz constant
of the product of bounded Lipschitz functions, and to the Leibniz rule
for differentiating the product of two functions.

        Suppose that $a_1$, $b_1$ are real numbers such that $a \le
a_1 \le b_1 \le b$.  If $f$ has bounded variation on $[a, b]$, then
$f$ has bounded variation on $[a_1, b_1]$, and
\begin{equation}
        V_{a_1}^{b_1}(f) \le V_a^b(f).
\end{equation}
This is because every partition of $[a_1, b_1]$ can be extended to a
partition of $[a, b]$.  In particular, $f$ is bounded on $[a, b]$ when
it has bounded variation.

        A partition $\mathcal{P}'$ of $[a, b]$ is said to be a
\emph{refinement} of a partition $\mathcal{P}$ of $[a, b]$ if
$\mathcal{P}'$ contains all of the terms in $\mathcal{P}$.
In this case, one can check that
\begin{equation}
        V_\mathcal{P}(f) \le V_{\mathcal{P}'}(f)
\end{equation}
for every $f : [a, b] \to {\bf R}$, using the triangle inequality.
Also, any finite collection of partitions of $[a, b]$ has a common
refinement.

        Suppose that $f$ has bounded variation on $[a, b]$, and that
$x \in (a, b)$.  Thus the restrictions of $f$ to $[a, x]$ and to $[x,
b]$ have bounded variation, and moreover
\begin{equation}
        V_a^x(f) + V_x^b(f) = V_a^b(f).
\end{equation}
Indeed, any partitions $\mathcal{P}_1$, $\mathcal{P}_2$ of $[a, x]$,
$[x, b]$, respectively, can be combined to get a partition
$\mathcal{P}_3$ of $[a, b]$ for which
\begin{equation}
        V_{\mathcal{P}_1}(f) + V_{\mathcal{P}_2}(f) = V_{\mathcal{P}_3}(f),
\end{equation}
which implies that $V_a^x(f) + V_x^b(f) \le V_a^b(f)$.  To get the
opposite inequality, note that every partition of $[a, b]$ can be
refined if necessary to contain $x$, and hence to be a combination of
partitions of $[a, x]$ and $[x, b]$.  The same argument shows that $f$
has bounded variation on $[a, b]$ if it has bounded variation on $[a,
x]$ and on $[x, b]$.

        Suppose that $f$ is continuously differentiable on $[a, b]$.
If $a \le r \le t \le b$, then
\begin{equation}
        |f(t) - f(r)| = \biggl|\int_r^t f'(\xi) \, d\xi\biggr|
                      \le \int_r^t |f'(\xi)| \, d\xi.
\end{equation}
This implies that
\begin{equation}
        V_\mathcal{P}(f) \le \int_a^b |f'(\xi)| \, d\xi
\end{equation}
for every partition $\mathcal{P}$ of $[a, b]$.  One can show that
\begin{equation}
        V_a^b(f) = \int_a^b |f'(\xi)| \, d\xi,
\end{equation}
using very fine partitions $\mathcal{P}$ of $[a, b]$.

        For each $r \in {\bf R}$, put $r_+ = r$ when $r \ge 0$ and
$r_+ = 0$ when $r \le 0$, and $r_- = -r$ when $r \le 0$ and $r-_ = 0$
when $r \ge 0$, so that
\begin{equation}
        r_+ - r_- = r, \quad r_+ + r_- = |r|.
\end{equation}
Given $f : [a, b] \to {\bf R}$ and a partition $\mathcal{P} =
\{t_j\}_{j = 0}^n$ of $[a, b]$, put
\begin{equation}
        P_\mathcal{P}(f) = \sum_{j = 1}^n (f(t_j) - f(t_{j - 1}))_+
\end{equation}
and
\begin{equation}
        N_\mathcal{P}(f) = \sum_{j = 1}^n (f(t_j) - f(t_{j - 1}))_-.
\end{equation}
Thus
\begin{equation}
        P_\mathcal{P}(f) + N_\mathcal{P}(f) = V_\mathcal{P}(f)
\end{equation}
and
\begin{equation}
        P_\mathcal{P}(f) - N_\mathcal{P}(f) = f(b) - f(a).
\end{equation}
Suppose that $f$ has bounded variation on $[a, b]$, and put
\begin{equation}
        P_a^b(f) = \sup \{P_\mathcal{P}(f) : \mathcal{P}
                                       \hbox{ is a partition of } [a, b]\}
\end{equation}
and
\begin{equation}
        N_a^b(f) = \sup \{N_\mathcal{P}(f) : \mathcal{P}
                                       \hbox{ is a partition of } [a, b]\}.
\end{equation}
One can check that
\begin{equation}
        P_a^b(f) + N_a^b(f) = V_a^b(f)
\end{equation}
and
\begin{equation}
        P_a^b(f) - N_a^b(f) = f(b) - f(a).
\end{equation}
Similarly,
\begin{equation}
        P_a^x(f) - N_a^x(f) = f(x) - f(a)
\end{equation}
when $a \le x \le b$.  This implies that $f$ can be expressed as the
difference of two monotone increasing functions on $[a, b]$, since
$P_a^x(f)$, $N_a^x(f)$ are monotone increasing in $x$.

        Functions of bounded variation do not have to be continuous,
but they can only have jump discontinuities.  More precisely, if $f$
has bounded variation on $[a, b]$, then $f$ has one-sided limits from
both sides at every point in $(a, b)$, and from the right and left
sides at $a$, $b$, respectively.  This follows from the analogous
statement for monotone functions and the fact that a function of
bounded variation can be expressed in terms of monotone functions, and
it can also be shown more directly.

\section{Lengths of paths}
\label{lengths of paths}
\setcounter{equation}{0}

        Let $(M, d(x, y))$ be a metric space, let $a$, $b$ be real
numbers with $a \le b$, and let $f$ be a function on $[a, b]$ with
values in $M$.  For each partition $\mathcal{P} = \{t_j\}_{j = 0}^n$
of $M$, consider
\begin{equation}
        \Lambda_\mathcal{P}(f) = \sum_{j = 1}^n d(f(t_j), f(t_{j - 1})).
\end{equation}
This is the same as the variation $V_\mathcal{P}(f)$ of $f$ on
$\mathcal{P}$ when $M$ is the real line with the standard metric.
If there is an upper bound for $\Lambda_\mathcal{P}$ over all
partitions $\mathcal{P}$ of $[a, b]$, then we say that the path $f :
[a, b] \to M$ has finite length, and the length of the path is defined by
\begin{equation}
         \Lambda_a^b(f) = \sup \{\Lambda_\mathcal{P} : \mathcal{P}
                                          \hbox{ is a partition of } [a, b]\}.
\end{equation}
This is the same as the total variation $V_a^b(f)$ of $f$ when $M =
{\bf R}$.  As in the previous case, $\Lambda_a^b(f) = 0$ if and only
if $f$ is constant.  If $a \le r \le t \le b$, then
\begin{equation}
        d(f(r), f(t)) \le \Lambda_\mathcal{P}(f)
\end{equation}
for any partition $\mathcal{P}$ of $[a, b]$ that contains $r$, $t$.
Hence $f([a, b])$ is a bounded set in $M$ when $f : [a, b] \to M$ has
finite length, with
\begin{equation}
        \diam f([a, b]) \le \Lambda_a^b(f).
\end{equation}
Of course, $f([a, b])$ is a compact set in $M$ when $f$ is continuous,
and therefore bounded.  If $f$ is continuous, then $f([a, b])$ is also
a connected set in $M$.

        If $f : [a, b] \to M$ is $C$-Lipschitz, then $f$ has finite
length, and
\begin{equation}
        \Lambda_a^b(f) \le C \, (b - a).
\end{equation}
Let $(\widetilde{M}, \widetilde{d}(u, v))$ be another metric space,
and suppose that $\phi : M \to \widetilde{M}$ is $C$-Lipschitz.  For
any $f : [a, b] \to M$ and partition $\mathcal{P}$ of $[a, b]$,
\begin{equation}
 \widetilde{\Lambda}_\mathcal{P}(\phi \circ f) \le C \, \Lambda_\mathcal{P}(f),
\end{equation}
where $\widetilde{\Lambda}$ is the analogous quantity for $\widetilde{M}$.
If $f : [a, b] \to M$ has finite length, then $\phi \circ f : [a, b]
\to \widetilde{M}$ does too, and
\begin{equation}
        \widetilde{\Lambda}_a^b(\phi \circ f) \le C \, \Lambda_a^b(f).
\end{equation}
In particular, if $f$ has finite length and $\phi : M \to {\bf R}$ is
Lipschitz, then $\phi \circ f$ has bounded variation.

        If $f : [a, b] \to M$ has finite length and $a \le a_1 \le b_1
\le b$, then the restriction of $f$ to $[a_1, b_1]$ has finite length, and
\begin{equation}
        \Lambda_{a_1}^{b_1}(f) \le \Lambda_a^b(f).
\end{equation}
If $\mathcal{P}$, $\mathcal{P}'$ are partitions of $[a, b]$ and
$\mathcal{P}'$ is a refinement of $\mathcal{P}$, then
\begin{equation}
        \Lambda_\mathcal{P}(f) \le \Lambda_{\mathcal{P}'}(f)
\end{equation}
for any $f : [a, b] \to M$, as in the case of real-valued functions in
the previous section.  As before, one can use this to show that
\begin{equation}
        \Lambda_a^x(f) + \Lambda_x^b(f) = \Lambda_a^b(f)
\end{equation}
for every $x \in (a, b)$ when $f$ has finite length.  If $a \le x <
b$, then
\begin{equation}
        \lim_{y \to x+} \Lambda_a^y(f)
\end{equation}
exists, because $\Lambda_a^x(f)$ is monotone increasing in $x$, and hence
\begin{equation}
        \lim_{y \to x+} \sup \{\Lambda_w^y(f) : x < w \le y\} = 0.
\end{equation}
This implies that
\begin{equation}
        \lim_{y \to x+} \diam f((x, y]) = 0,
\end{equation}
since
\begin{eqnarray}
        \diam f((x, y]) & = & \sup \{\diam f([w, y]) : x < w \le y\}  \\
             & \le & \sup \{\Lambda_w^y(f) : x < w \le y\}. \nonumber
\end{eqnarray}
If $M$ is complete, then it follows that $f$ has a limit from the
right at $x$, and similarly there is a limit from the left when $a < x
\le b$.

        Suppose now that $M$ is ${\bf R}^n$, equipped with a norm $N$,
and thus the metric $d_N(x, y)$ associated to $N$ too.  If $f_1, f_2 :
[a, b] \to {\bf R}^n$ have finite length and $r_1, r_2 \in {\bf R}$,
then $r_1 \, f_1 + r_2 \, f_2$ has finite length, and
\begin{equation}
        \Lambda_a^b(r_1 \, f_1 + r_2 \, f_2)
                   \le |r_1| \, \Lambda_a^b(f_1) + |r_2| \, \Lambda_a^b(f_2).
\end{equation}
This is similar to the case of real-valued functions, and one can also
treat the product of a real-valued function and an ${\bf R}^n$-valued
function on $[a, b]$ in the same way as before.  If $f : [a, b] \to
{\bf R}^n$ is continuously differentiable, then one can show that $f$
has finite length and that
\begin{equation}
         \Lambda_a^b(f) = \int_a^b N(f'(\xi)) \, d\xi,
\end{equation}
in practically the same way as before.  It can be interesting to
consider integral norms
\begin{equation}
         \Big(\int_a^b N(f'(\xi))^p \, d\xi\Big)^{1/p}
\end{equation}
as well, $1 \le p < \infty$.  The $p = \infty$ case corresponds to the
maximum of $N(f')$ on $[a, b]$.  This integral norm is especially
interesting when $p = 2$ and $N$ is the standard Euclidean norm on
${\bf R}^n$.  For other $p$, there is some simplification when $N(v) =
\|v\|_p$.  If $N(v) = \|v\|_1$, then the length of any path of finite
length in ${\bf R}^n$ is equal to the sum of the total variations of
the coordinates of the path.  This uses the fact that any finite
collection of partitions of $[a, b]$ has a common refinement, so that
independent partitions for the coordinate functions are equivalent to
using the same partition for the whole path.

\section{Snowflake metrics}
\label{snowflakes}
\setcounter{equation}{0}

        Let $\alpha$ be a positive real number, with $\alpha < 1$.
For any pair of nonnegative real numbers $u$, $v$,
\begin{equation}
\label{u^alpha, v^alpha}
        (u + v)^\alpha \le u^\alpha + v^\alpha.
\end{equation}
To see this, observe that
\begin{equation}
        \max(u, v) \le (u^\alpha + v^\alpha)^{1/\alpha},
\end{equation}
and hence
\begin{equation}
        u + v \le \max(u, v)^{1 - \alpha} \, (u^\alpha + v^\alpha)
               \le (u^\alpha + v^\alpha)^{1/\alpha}.
\end{equation}
Note that the inequality is strict in (\ref{u^alpha, v^alpha}) when $u, v > 0$.

        If $(M, d(x, y))$ is a metric space, then it follows from
(\ref{u^alpha, v^alpha}) that $d(x, y)^\alpha$ is also a metric on
$M$.  This does not change the topology of $M$, but it does change the
geometry.  Many familiar examples of snowflake curves in the plane
have approximately this type of geometry, for instance.

        Suppose that $f : [a, b] \to M$ is a continuous path with
finite length with respect to $d(x, y)^\alpha$.  This means that
\begin{equation}
        \sum_{j = 1}^n d(f(t_j), f(t_{j - 1}))^\alpha \le A
\end{equation}
for some $A \ge 0$ and every partition $\{t_j\}_{j = 0}^n$ of $[a, b]$.
Let $\epsilon > 0$ be given.  By continuity and compactness, $f$
is uniformly continuous, and so there is a $\delta > 0$ such that
\begin{equation}
        d(f(r), f(w)) < \epsilon
\end{equation}
for every $r, w \in [a, b]$ such that $|r - w| < \delta$.  Hence
\begin{equation}
        \sum_{j = 1}^n d(f(t_j), f(t_{j - 1})) \le \epsilon^{1 - \alpha} \, A
\end{equation}
when $t_j - t_{j - 1} < \delta$ for each $j = 1, \ldots, n$.  Every
partition of $[a, b]$ has a refinement with this property, which
implies that the length $\Lambda_a^b(f)$ of $f$ with respect to $d(x,
y)$ satisfies
\begin{equation}
        \Lambda_a^b(f) \le \epsilon^{1 - \alpha} \, A.
\end{equation}
Thus $\Lambda_a^b(f) = 0$, since $\epsilon > 0$ is arbitrary, and $f$
must be constant.

\section{H\"older continuity}
\label{holder}
\setcounter{equation}{0}

        Let $(M_1, d_1(x, y))$ and $(M_2, d_2(w, z))$ be metric
spaces.  A mapping $f : M_1 \to M_2$ is said to be \emph{H\"older
continuous} of order $\alpha$, $0 < \alpha < 1$, if
\begin{equation}
        d_2(f(x), f(y)) \le C \, d_1(x, y)^\alpha
\end{equation}
for some $C \ge 0$ and every $x, y \in M_1$.  One might also say that
$f$ is Lipschitz of order $\alpha$ in this case, but it will be
convenient to refer to this as a Lipschitz condition when $\alpha = 1$
and H\"older continuity when $0 < \alpha < 1$.  Similar names are
sometimes used for other related conditions as well.

        As in the previous section, $d_1(x, y)^\alpha$ is a metric on
$M_1$, and therefore $f$ is H\"older continuous of order $\alpha$ with
respect to $d_1(x, y)$ if and only if $f$ is Lipschitz with respect to
$d_1(x, y)^\alpha$.  Thus many basic properties of H\"older continuous
mappings follow from the corresponding statements for Lipschitz
mappings.  In particular,
\begin{equation}
        f_p(x) = d_1(p, x)^\alpha
\end{equation}
is a real-valued H\"older continuous function of order $\alpha$ on
$M_1$ with $C = 1$ for each $p \in M_1$.

        Let $(M, d(x, y))$ be a metric space, and consider the case of
a continuous path $f : [a, b] \to M$.  If $f$ is $C$-Lipschitz, then
$f$ is H\"older continuous of order $\alpha$ for each $\alpha \in (0,
1)$ with constant $C \, (b - a)^{1 - \alpha}$.  Of course, $f$ also
has finite length $\le C \, (b - a)$ when $f$ is $C$-Lipschitz.
However, continuous paths of finite length need not be H\"older
continuous of any positive order, and there are couterexamples already
for monotone increasing real-valued functions.  Similarly, H\"older
continuous paths may not have finite length.

        Suppose that $f : [a, b] \to M$ is H\"older continuous of
order $\alpha$ with constant $C > 0$.  For each $\rho > 0$, $f([a,
b])$ is contained in the union of $O(\rho^{-1/\alpha})$ subsets of $M$
with diameter $\le \rho$, because $[a, b]$ is the union of
$O(\rho^{-1/\alpha})$ subintervals of length $\le
(\rho/C)^{1/\alpha}$.  This implies that the Minkowski dimension of
$f([a, b])$ is $\le 1/\alpha$, and hence the Hausdorff dimension of
$f([a, b])$ is $\le 1/\alpha$ too.  This is an analogue for H\"older
continuous paths of the fact that Lipschitz paths have finite length.

\section{Coverings}
\label{coverings}
\setcounter{equation}{0}

        If $[a, b], [a_1, b_1], \ldots, [a_n, b_n]$ are closed
intervals in the real line such that
\begin{equation}
        [a, b] \subseteq \bigcup_{j = 1}^n [a_j, b_j],
\end{equation}
then
\begin{equation}
        b - a \le \sum_{j = 1}^n (b_j - a_j).
\end{equation}
More generally, if $E_1, \ldots, E_n$ are bounded subsets of ${\bf R}$
such that
\begin{equation}
        [a, b] \subseteq \bigcup_{j = 1}^n E_j,
\end{equation}
then
\begin{equation}
        b - a \le \sum_{j = 1}^n \diam E_j.
\end{equation}
Indeed, each $E_j$ is contained in a closed interval of the same diameter.

        Let $(M, d(x, y))$ be a metric space, and let $A, E_1, \ldots, E_n$
be bounded subsets of $M$ such that
\begin{equation}
        A \subseteq \bigcup_{j = 1}^n E_j.
\end{equation}
If $A$ is connected, then
\begin{equation}
        \diam A \le \sum_{j = 1}^n \diam E_j.
\end{equation}
To see this, remember first that continuous mappings send connected
sets to connected sets.  If $f : M \to {\bf R}$ is continuous, then
$f(A)$ is an interval in the real line, which may be open, or closed,
or half-open and half-closed.  At any rate,
\begin{equation}
        {\diam}_{\bf R} f(A) \le \sum_{j = 1}^n {\diam}_{\bf R} f(E_j),
\end{equation}
where the subscripts indicate that these are diameters in ${\bf R}$.
If $f$ is $C$-Lipschitz, then
\begin{equation}
        {\diam}_{\bf R} f(A) \le C \, \sum_{j = 1}^n \diam E_j.
\end{equation}
The desired estimate follows by applying this to $1$-Lipschitz
functions of the form $f_p(x) = d(p, x)$, $p \in A$.

        The hypothesis that $A$ be connected is essential here.  If
$A$ is a finite set with at least two elements, then $\diam A > 0$,
but $A$ is contained in the union of finitely many sets with one
element and thus diameter $0$.  Cantor's middle-thirds set in the real
line has diameter equal to $1$, and is contained in the union of
$2^\ell$ intervals of length $3^{-\ell}$ for each $\ell \ge 1$.  A
compact set $A \subseteq {\bf R}$ has Lebesgue measure $0$ exactly if
for each $\epsilon > 0$ there are finitely many bounded sets $E_1,
\ldots, E_n \subseteq {\bf R}$ such that $A \subseteq \bigcup_{j =
1}^n E_j$ and $\sum_{j = 1}^n \diam E_j < \epsilon$.

        In particular, if $A \subseteq M$ is a bounded connected set
and $\rho > 0$, then $A$ is not covered by fewer than $\diam A / \rho$
bounded subsets of $M$ of diameter $\le \rho$.  Depending on the
situation, many more of these subsets may be required.  If $M$ is
${\bf R}^n$ equipped with the standard metric, for example, then a
bounded set $A$ can be covered by $O(\rho^{-n})$ sets of diameter $\le
\rho$.  One needs at least a positive multiple of $\rho^{-n}$ such
sets when $A$ has nonempty interior, because otherwise the
$n$-dimensional volume of $A$ would be too small.

\section{Domains in ${\bf R}^n$}
\label{domains, R^n}
\setcounter{equation}{0}

        A set $U$ in a metric space $M$ is an open set if for every $p
\in U$ there is an $r > 0$ such that $B(p, r) \subseteq U$.  Any norm
$N$ on ${\bf R}^n$ determines the same open sets as the standard
metric.  This is because $N$ is less than or equal to a constant times
the standard norm, and conversely the standard norm is less than or
equal to a constant times $N$.  The first statement can be checked
directly by expressing any element of ${\bf R}^n$ as a linear
combination of the standard basis for ${\bf R}^n$ and using the
triangle inequality.  As mentioned previously, this and the triangle
inequality imply that $N$ is continuous with respect to the standard
norm.  Hence the minimum of $N$ is attained on the standard unit
sphere in ${\bf R}^n$, since the latter is compact.  The standard norm
times the minimum of $N$ on the unit sphere is less than or equal to
$N$ on all of ${\bf R}^n$, by homogeneity, which implies the second
statement.  For explicit norms like $\|w\|_p$, $1 \le p \le \infty$,
the comparison with the standard norm can be verified directly.

        Suppose that $U$ is a connected open set in ${\bf R}^n$, which
is to say that $U$ is not the union of two disjoint nonempty open
sets.  It is well known that $U$ is then pathwise-connected, so that
for every $p, q \in U$ there is a continuous mapping $f : [a, b] \to
U$ such that $f(a) = p$ and $f(b) = q$.  More precisely, one can even
take $f$ to be piecewise-affine on $[a, b]$.  In particular, $f$ then
has finite length.

        However, it is not clear how small the length of $f$ can be.
Of course, the length of $f$ is at least the distance between $p$ and
$q$.  If $U$ is convex, then one can take $f$ to be affine, and the
length of $f$ is equal to the distance between $p$ and $q$.
Otherwise, the length of $f$ may have to be quite large compared to
the distance between $p$ and $q$.  It is easy to give examples where
this happens in the plane.  For instance, there may be elements of $U$
on opposite sides of the boundary locally.  The boundary of $U$ might
also be complicated, with spirals or other obstacles.

        Even if the boundary of $U$ is complicated, it may be that $U$
behaves well in terms of lengths of paths.  For example, if $U$ is the
region in the plane bounded by the von Koch snowflake, then every pair
of elements of $U$ can be connected by a path of length bounded by a
constant multiple of the distance between them.  The main idea is for
the path to avoid the boundary as much as possible, without going too
far away.  There can also be relatively small parts of the boundary
that only cause minor detours for paths in the domain.

\section{Lipschitz graphs}
\label{lipschitz graphs}
\setcounter{equation}{0}

        Let $k, l, n$ be positive integers such that $k + l = n$, and
let us identify ${\bf R}^n$ with ${\bf R}^k \times {\bf R}^l$, so that
an element $x$ of ${\bf R}^n$ may be expressed as $(x', x'')$, where
$x' \in {\bf R}^k$ and $x'' \in {\bf R}^l$.  Also let $A : {\bf R}^k
\to {\bf R}^l$ be a continuous mapping, and consider its graph
\begin{equation}
        \{(x', x'') \in {\bf R}^n : x'' = A(x')\}.
\end{equation}
This is a nice $k$-dimensional topological submanifold of ${\bf R}^n$.
If $k = n - 1$, then this hypersurface has two complementary
components $U_+$, $U_-$ consisting of $(x', x'') \in {\bf R}^n$ such
that $x'' > A(x')$ and $x'' < A(x')$, respectively.  If $k < n$, then
the complement of the graph in ${\bf R}^n$ is connected.  For any $k$,
\begin{equation}
        (x', x'') \mapsto (x', x'' + A(x'))
\end{equation}
defines a homeomorphism on ${\bf R}^n$ that sends the $k$-plane $x'' =
0$ to the graph of $A$.  If $f(t)$ is a continuous path in ${\bf R}^k$
parameterized by an interval $[a, b]$, then $\widehat{f}(t) = (f(t),
A(f(t))$ is a continuous path in the graph of $A$.  The graph of $A$
is itself a curve in ${\bf R}^n$ when $k = 1$.

        Suppose that $A$ is Lipschitz.  If $f$ has finite length, then
$\widehat{f}$ does too, and the length of $\widehat{f}$ is bounded by
a constant multiple of the length of $f$.  If the Lipschitz constant
of $A$ is small, then this constant multiple is close to $1$.
Using affine paths in ${\bf R}^k$, we get that every pair of elements
of the graph of $A$ can be connected by a continuous path in the graph
of $A$ of finite length bounded by a constant multiple of the distance
between them, where the constant multiple is close to $1$ when $A$
has small Lipschitz constant.

        A $k$-dimensional $C^1$ submanifold of ${\bf R}^n$ is locally
the same as the graph of a continuously-differentiable mapping on
${\bf R}^k$ with respect to a suitable choice of coordinate axes.  By
rotating the axes so that ${\bf R}^k$ is parallel to the tangent plane
of the submanifold at a particular point, the submanifold can be
represented near the point as the graph of a function with small
Lipschitz constant.  The Lipschitz constant tends to $0$ as one
approaches the point in question.  Thus distances on $C^1$
submanifolds are approximately the same as the infimum of lengths of
paths on the submanifold locally.

\section{Real analysis}
\label{real analysis}
\setcounter{equation}{0}

        For the sake of simplicity, we have so far avoided referring
to Lebesgue integrals and measure.  However, this more sophisticated
theory can be quite convenient in the present context.  Let us mention
some of the key points.

        A basic fact is that a monotone real-valued function on an
open interval in the real line is differentiable ``almost
everywhere'', which is to say on the complement of a set of Lebesgue
measure $0$.  Thus additional hypotheses of differentiability are
sometimes superfluous.  Unfortunately, even continuous monotone
functions cannot necessarily be recoved from their almost everywhere
derivative, as in the fundamental theorem of calculus, without an
extra condition of ``absolute continuity''.  Indeed, there are
examples of nonconstant continuous monotone increasing functions with
derivative equal to $0$ almost everywhere.

        It follows that a real-valued function of bounded variation on
an interval in the real line is differentiable almost everywhere,
since it can be expressed as the difference of two monotone increasing
functions.  In particular, a real-valued Lipschitz function on an
interval is differentiable almost everywhere.  Lipschitz functions are
absolutely continuous, and so there is a version of the fundamental
theorem of calculus for them.  As corollaries of this fact, a
Lipschitz function $f$ on an interval is constant if $f'(x) = 0$
almost everywhere, $f$ is monotone increasing if $f'(x) \ge 0$ almost
everywhere, and $f$ is $C$-Lipschitz if $|f'(x)| \le C$ almost
everywhere.

        At the same time, bounded variation and Lipschitz conditions
have natural extensions involving metric spaces, as we have seen.  The
composition of a path of finite length in a metric space with a
real-valued Lipschitz function on the metric space is a function of
bounded variation, which is Lipschitz when the path is.  There are
also a lot of real-valued Lipschitz functions on any metric space.
Even on ${\bf R}^n$, there are a lot of nice functions that are
Lipschitz and not continuously differentiable, such as the distance to
a point or to a set.

\end{document}